\documentclass[11pt,leqno]{article}
\usepackage{amsthm,amsfonts,amssymb,amsmath,oldgerm}
\usepackage{epsfig}
\numberwithin{equation}{section} 

\def\uref{{u^{\rm ref}}}
\def\alpharef{{\alpha^{\rm ref}}}

\def\eps{\varepsilon }

\newcommand\R{\mathbb R}

\def\eps{\varepsilon}

\newcommand\errfn{\textrm{errfn}}
\newcommand\br{\begin{remark}}
\newcommand\er{\end{remark}}
\newcommand\bp{\begin{pmatrix}}
\newcommand\ep{\end{pmatrix}}
\newcommand\be{\begin{equation}}
\newcommand\ee{\end{equation}}
\newcommand\ba{\begin{equation}\begin{aligned}}
\newcommand\ea{\end{aligned}\end{equation}}

\newcommand{\bap}{\begin{app}}
\newcommand{\eap}{\end{app}}
\newcommand{\begs}{\begin{exams}}
\newcommand{\eegs}{\end{exams}}
\newcommand{\beg}{\begin{example}}
\newcommand{\eeg}{\end{exaplem}}
\newcommand{\bpr}{\begin{proposition}}
\newcommand{\epr}{\end{proposition}}
\newcommand{\bt}{\begin{theorem}}
\newcommand{\et}{\end{theorem}}
\newcommand{\bc}{\begin{corollary}}
\newcommand{\ec}{\end{corollary}}
\newcommand{\bl}{\begin{lemma}}
\newcommand{\el}{\end{lemma}}
\newcommand{\bd}{\begin{definition}}
\newcommand{\ed}{\end{definition}}
\newcommand{\brs}{\begin{remarks}}
\newcommand{\ers}{\end{remarks}}

\newtheorem{theo}{Theorem}[section]

\newtheorem{exams}[theo]{Examples}

\numberwithin{equation}{section}

\newcommand{\RR}{{\mathbb R}}

\newtheorem{theorem}{Theorem}[section]
\newtheorem{proposition}[theorem]{Proposition}
\newtheorem{corollary}[theorem]{Corollary}
\newtheorem{lemma}[theorem]{Lemma}
\newtheorem{definition}[theorem]{Definition}

\newtheorem{example}[theorem]{Example}
\newtheorem{remark}[theorem]{Remark}

\newcommand\cK{{\cal  K}}

\title{
Instantaneous shock location and one-dimensional nonlinear stability of
viscous shock waves} 

\author{\sc \small 
Kevin Zumbrun
\thanks{Indiana University, Bloomington, IN 47405;
kzumbrun@indiana.edu:
Research of K.Z. was partially supported
under NSF grants no. DMS-0300487 and DMS-0801745.
 }
}
\begin{document}

\maketitle


\begin{abstract}
We illustrate in a simple setting
the instantaneous shock tracking approach 
to stability of viscous conservation laws
introduced by Howard, Mascia, and Zumbrun.
This involves a choice of the definition
of instanteous location of a viscous shock--
we show that this choice is time-asymptotically
equivalent both to the natural choice of least-squares
fit pointed out by Goodman and to a simple phase condition used by
Gu\`es, M\'etivier, Williams, and Zumbrun in
other contexts.
More generally, we show that it is asymptotically
equivalent to any location defined by a localized projection.
\end{abstract}


\bigbreak
\section{Introduction}\label{s:intro}
In this note,
we illustrate in the simple and concrete setting of Burgers equation the
argument for nonlinear stability of viscous shock waves developed for
general systems of conservation laws in \cite{Z1,MaZ2,MaZ3,MaZ4},
based on instantaneous tracking of the location of the perturbed viscous
shock wave.
The advantage of Burgers equation is that the linearized equations
may be solved explicitly by a linearized Hopf--Cole transformation,
thus isolating the nonlinear issues we wish to discuss.
This same example was given in \cite{Z1}; here we 
expand a bit the
surrounding discussion, reexamining the question of what is 
a reasonable or natural definition of the instantaneous location
of a perturbed viscous shock wave and adding a discussion
of the small-amplitude limit.

Using the purely operational but analytically tractable definition of 
\cite{Z1,MaZ2,MaZ4} as a tool
for comparison, we show that {\it any} definition based on  localized projection
is time-asymptotically equivalent to any other and to the 
definition of \cite{Z1,MaZ2,MaZ4}; 
see Appendix \ref{s:alt}, and especially 
Remarks \ref{trackrmk1}--\ref{trackrmk2}.
Moreover, any of these may be used as the basis of a nonlinear stability 
argument.
This includes in particular both the natural definition by least squares fit pointed out early on by Goodman \cite{G}
and, in the limit of infinite localization (to a single point), a very simple definition
based on a phase condition, introduced by Gu\`es, M\'etivier, Williams, and Zumbrun \cite{GMWZ}.

Our analysis is {\it intended for the nonspecialist}.
It is brief and self-contained except for standard 
linear and short-time parabolic existence theory.
Though we restrict here for simplicity to the scalar Burgers case, our arguments and conclusions 
extend in straightforward fashion to the general system case
\cite{Z1,MaZ2,MaZ4}, {\it once there are established
the requisite bounds on the linearized solution operator}.
This separate, and in general difficult, problem has been treated in 
\cite{ZH,MaZ3,Z2}; see Remark \ref{stabrmk} and the discussion
of Section \ref{s:sys}.
Our purpose here is, rather, to isolate the issues connected with 
viscous shock-tracking and the nonlinear iteration argument
by restricting to a case where the linearized bounds are
available by exact solution formula.

\subsection{Problem and equations}
Consider the scalar viscous conservation law
\begin{equation}\label{burgers}
u_t + f(u)_x=u_{xx},
\end{equation}
$u=u(x,t)\in \R$, $x\in \R$, $t\in \R^+$, with
\begin{equation}\label{f}
f(u)=u^2/2.
\end{equation}
Eq. \eqref{burgers} serves as a simple model for gas dynamics, traffic flow, or
shallow-water waves, where $u$ represents the density of some
conserved quantity and $f$ its flux through a fixed point $x$.
With the choice of flux \eqref{f}, \eqref{burgers} becomes {\it Burgers equation},
the prototypical example of a scalar viscous conservation law.
Behavior for other (convex) fluxes is qualitatively similar.

We investigate the question of
{\it nonlinear stability} of solutions $u=\bar u$, that is, 
whether a perturbation $\tilde u=\bar u+u$ remains close to $\bar u$ in some norm
for initial perturbations $u|_{t=0}=(\tilde u-\bar u)|_{t=0}$ sufficiently small
in some (possibly different) norm:
more specifically, {\it nonlinear asymptotic stability}, that is,
whether $\tilde u$ both remains near and converges to $\bar u$ as $t\to +\infty$ for
initial perturbations sufficiently small.
Since the equation \eqref{burgers} is translation-invariant, we must 
when relevant (specifically, when translates of $\bar u$ are
not equal to $\bar u$) 
adjust the second notion
to that of {\it nonlinear asymptotic orbital stability}, 
defined as nonlinear stability together with convergence as $t\to +\infty$
to the set of translates of $\bar u$,
as discussed further below.

\subsection{Constant and traveling-wave solutions}
An obvious class of solutions of \eqref{burgers} are
the set of constant solutions $\bar u\equiv a$, $a\in \RR$.
A second class of solutions are {\it viscous shock waves},
or smooth traveling-wave solutions
\begin{equation}\label{op}
u(x,t)=\bar u(x-st),
\qquad
\lim_{x\to \pm \infty}\bar u(x)=u_\pm
\end{equation}
$s$ constant, connecting constant endstates $u_\pm$.
If $s=0$, they are equilibria, or {\it stationary waves}
of the associated evolution equation \eqref{burgers}.
A traveling wave may always be converted to a standing wave by
the change of coordinates
$x\to x-st$ to  a frame moving with the same speed $s$.


Observing that 
$ \partial_t \bar u(x-st)=-s\bar u'$,
$\partial_x \bar u(x-st)=\bar u'$,
and $\partial_x^2 \bar u(x-st)=\bar u'',$
we obtain for a solution \eqref{op} 
the {\it profile equation} $ -s\bar u'+ f(\bar u)'=\bar u''.  $
Integrating from $-\infty$ to $x$ reduces this to a first-order equation
\be\label{prof}
\bar u'= (f(\bar u)-s\bar u)-
(f(u_-)-su_-).
\ee
For definiteness taking $s=0$, $u_-=1$, we obtain 
$
\bar u'= (1/2)(1-\bar u^2),
$
which has the explicit solution
\begin{equation}\label{tanhform}
\bar u(x)=-\tanh(x/2)
\end{equation}
connecting endstates $u_\pm=\mp 1$.
This is the unique solution up to translation in $x$
connecting that particular pair of endstates.
Other endstates and speeds also lead to $\tanh$ profiles,
as may be seen by invariances of Burgers equation; thus,
we may without loss of generality restrict to this specific case.

\section{Stability of constant solutions}\label{s:const}
To indicate the basic approach,
let us first consider stability of a constant solution 
\be\label{asoln}
 \bar u\equiv a,\qquad a\in \RR
\ee
of \eqref{burgers}.  
Letting $\tilde u$ be a second solution of \eqref{burgers}, and
defining perturbation $u:=\tilde u-\bar u$,
we obtain after a brief computation the {\it perturbation equation}
\be\label{cpert}
u_t-Lu= N(u)_x,
\ee
where $Lu:= u_{xx} -au_x $
is the linearization of $u_{xx}-f(u)_x$ about solution $\bar u\equiv a$.
and $N(u):= -u^2/2$ is a quadratic order remainder.


\subsection{Linear solution operator}

The homogeneous linearized equations $v_t-Lv=0$ may be recognized
as a {\it convected heat equation}
\be\label{heat}
v_t+av_x=v_{xx},
\qquad v|_{t=0}=f .
\ee
This admits an exact solution
\be\label{linsol}
e^{Lt}f=\int_{-\infty}^{+\infty} G(x,t;y)f(y)dy,
\ee
where 
\be\label{lGreen}
G(x,t;y):=e^{Lt}\delta_y(x)=
\frac{ e^{-\frac{|x-y-at|^2}{4t}}}
{\sqrt{4 \pi t} }
\ee
is the Green function for \eqref{heat}, a convected heat-kernel.
This yields in particular a unique classical solution 
$v\in C^0(t\ge 0;L^p(x)) \cap C^2(t>0,x)$ 
for each $f\in L^p$.

Easy scaling arguments yields, for $1\le p\le \infty$,
\ba\label{scale}
|G(\cdot, t;y)|_{L^p(x)}&=
|G(x, t;\cdot)|_{L^p(y)}= C_pt^{-\frac{1}{2}(1-1/p)},\\
|G_y(\cdot, t;y)|_{L^p(x)}&=
|G_y(x, t;\cdot)|_{L^p(y)}= C_p't^{-\frac{1}{2}(1-1/p)-\frac{1}{2}},\\
\ea
for some constants $C_p$, $C_p'>0$.
From \eqref{scale}, we readily obtain the
following {\it linearized estimates}
(standard heat kernel bounds).
\bl\label{l:heat}
For some $C>0$, all $t>0$,
\ba\label{heatbds}
\Big|\int_{-\infty}^{+\infty} G(x,t;y)f(y)dy\Big|_{L^p}&\le 
Ct^{-\frac{1}{2}(1-1/p)}|f|_{L^1(x)},\\
\Big|\int_{-\infty}^{+\infty} G_y(x,t;y)f(y)dy\Big|_{L^p}&\le 
Ct^{-\frac{1}{2}(1-1/p)-\frac{1}{2}}|f|_{L^1(x)}.\\
\ea
\el

\begin{proof}
Applying the Triangle inequality together with \eqref{scale}, we obtain
$$
\Big|\int_{-\infty}^{+\infty} G(x,t;y)f(y)dy\Big|_{L^p(x)}\le 
\int_{-\infty}^{+\infty} |G(\cdot,t;y|_{L^p(x)}|f(y)|dy
= C_pt^{-\frac{1}{2}(1-1/p)}|f|_{L^1}.
$$
The proof of the second inequality is similar.
\end{proof}

\bl\label{l:heat2}
For some $C>0$, all $t>0$,
\ba\label{heatbds2}
\Big|\int_{-\infty}^{+\infty} G(x,t;y)f(y)dy\Big|_{L^p(x)}&\le 
C |f|_{L^p},
\quad
\Big|\int_{-\infty}^{+\infty} G_y(x,t;y)f(y)dy\Big|_{L^p(x)}\le 
Ct^{-\frac{1}{2}}|f|_{L^p}.\\
\ea
\el

\begin{proof}
Noting that $G(x,t;y)=G(x-y,t;0)$, so that
\eqref{linsol} is a convolution,
we may rewrite 
$\int_{-\infty}^{+\infty} G(x,t;y)f(y)dy$ as
$\int_{-\infty}^{+\infty} G(z,t;0)f(x-z)dz$ with $z:=x-y$. 
Applying the Triangle inequality and \eqref{scale}, we obtain
$$
\Big|\int_{-\infty}^{+\infty} G(x,t;y)f(y)dy\Big|_{L^p(x)}\le 
\int_{-\infty}^{+\infty} |G(z,t;0||f|_{L^p}dz
= C_1 |f|_{L^p}.
$$
The proof of the second inequality is similar.
\end{proof}

\subsection{Integral representation}
From the homogeneous linearized solution formula
\eqref{linsol}, we obtain by variation of constants/Duhamel's formula
a solution for the inhomogeneous linearized equations 
\be\label{linhom}
v_t-Lv=g,\quad v|_{t=0}=f
\ee
of $v=e^{Lt}f +\int_0^t e^{L(t-s)}g(s)ds$, or
\be\label{ilinsol}
v(x,t)=
\int_{-\infty}^{+\infty} G(x,t;y)f(y)dy +
\int_0^t \int_{-\infty}^{+\infty} G(x,t-s;y)g(y,s)dy\, ds,
\ee
yielding a unique $C^0(t\ge 0;L^p(x))\cap C^2(t>0;x)$ solution $v$
for $f\in L^p$ and $g\in W^{-1,p}$.

\subsection{Nonlinear iteration}
Returning now to the nonlinear problem \eqref{cpert}, we
have, setting $g=N(u)_x$ in \eqref{linhom}, the representation
$u(x,t)=
\int_{-\infty}^{+\infty} G(x,t;y)u_o(y)dy +
\int_0^t \int_{-\infty}^{+\infty} G(x,t-s;y)N(u(y,s))_y dy\, ds$,
or, integrating the last term by parts,
\be\label{nintrep}
u(x,t)=
\int_{-\infty}^{+\infty} G(x,t;y)u_0(y)dy 
-\int_0^t \int_{-\infty}^{+\infty} G_y(x,t-s;y)N(u(y,s)) dy\, ds,
\ee
valid so long as the solution $u$ exists and remains sufficiently smooth
that \eqref{nintrep} gives the unique solution to the associated
inhomogeneous problem, in particular for $u_0\in L^p\cap L^\infty$
and $u$ in $C^0(t\ge 0; L^p\cap L^\infty)\cap C^2(t>0,x)$,
any $p\ge 1$.

On the other hand, standard short-time existence theory (proved, e.g.,
by contraction-mapping using a similar representation with shifted
initial time) yields existence of a $C^0(t\ge 0;L^p\cap L^\infty(x))
\cap C^2(t>0,x)$ solution so long
as $|u|_{L^p\cap L^\infty}$ remains bounded.

Define now
\be\label{czeta}
\zeta(t):=\sup_{0\le s\le t,\, 1\le p\le \infty} |u|_{L^p}(s)(1+t)^{\frac{1}{2}(1-1/p)}.
\ee

\bl\label{cclaim}
For all $t\ge 0$ for which $\zeta(t)$ is finite, some $C>0$,
and $E_0:=|u_0|_{L^1\cap L^\infty}$,
\be\label{eq:cclaim}
\zeta(t)\le C(E_0+\zeta(t)^2).
\ee
\el

\begin{proof}
Noting, by quadratic dependence $N(u)=O(|u|^2)$ and the definition
\eqref{czeta} of $\zeta$, that 
\ba\label{Nbds}
|N(u)|_{L^1}&\le C|u|_{L^2}^2\le \zeta(t)^2(1+t)^{-\frac{1}{2}}\\
|N(u)|_{L^p}&\le C|u|_{L^p}|u|_{L^\infty} \le 
\zeta(t)^2(1+t)^{-\frac{1}{2}(1-1/p)-\frac{1}{2}} ,
\ea
we obtain, applying Lemmas \ref{heatbds}--\ref{heatbds2} to representation 
\eqref{nintrep}, the estimate
\ba\label{cest}
|u(\cdot,t)|_{L^p(x)}& \le
\Big|\int_{-\infty}^{+\infty} G(x,t;y)u_0(y)dy \Big|_{L^P(x)}\\
&\quad
+\Big|\int_0^{t/2} \int_{-\infty}^{+\infty} 
G_y(x,t-s;y)N(u(y,s)) dy\, ds\Big|_{L^P(x)}\\
&\quad +\Big|\int_{t/2}^t \int_{-\infty}^{+\infty} 
G_y(x,t-s;y)N(u(y,s)) dy\, ds\Big|_{L^P(x)}\\
&\le
C(1+t)^{-\frac{1}{2}(1-1/p)}E_0 +
C\zeta(t)^2\int_0^{t/2} (t-s)^{-\frac{1}{2}(1-1/p)-1/2}
(1+s)^{-\frac{1}{2}}ds\\
&\quad
+C\zeta(t)^2\int_{t/2}^{t} (t-s)^{-\frac{1}{2}}
(1+s)^{-\frac{1}{2}(1-1/p)-\frac{1}{2}}ds\\
&
\le
 C(E_0+\zeta(t)^2) (1+t)^{-\frac{1}{2}(1-1/p)}.
\ea
Rearranging, we obtain \eqref{cclaim}.
\end{proof}

\bc[Stability of constant solutions]\label{ccor}
Constant solutions $\bar u\equiv a$ are nonlinearly 
stable in $L^1\cap L^\infty$ and nonlinearly asympotically
stable in $L^p$, $p>1 $, with respect to 
initial perturbations $u_0$ that are sufficiently small in $L^1\cap L^\infty$.
More precisely, for some $C>0$,
\be\label{eq:cest}
|\tilde u-\bar u|_{L^p}(t)\le C(1+t)^{-\frac{1}{2}(1-1/p)}
|\tilde u-\bar u|_{L^1\cap L^\infty}|_{t=0}
\ee
for all $t\ge 0$, $1\le p\le \infty$, 
for solutions $\tilde u$ of \eqref{burgers} with
$|\tilde u-\bar u|_{L^1\cap L^\infty}|_{t=0}$ sufficiently small.
\ec

\begin{proof}
(``Continuous induction'')
By Lemma \ref{cclaim}, $\zeta(t)\le C(E_0+\zeta(t)^2)$
for 
\be\label{Edef}
 E_0:= |\tilde u-\bar u|_{L^1\cap L^\infty}|_{t=0}.
\ee
Taking $E_0< \frac{1}{4C^2}$, we have therefore that
$\zeta(t)< 2CE_0$ whenever $\zeta(t) \le 2CE_0$,
and so the set of $t\ge 0$ for which $\zeta(t)< 2CE_0$ is
equal to the set of $t\ge 0$ for which $\zeta(t)\le 2CE_0$.
Recalling,
by the cited standard short-time existence theory,
that $\zeta$ is continous wherever it is finite,
we find, therefore, that
the set of $t\ge 0$ for which $\zeta(t)< 2CE_0$ is both open 
and closed.
Taking without loss of generality $C>1/2$, so that $t=0$ is
contained in this set, we have that the set is nonempty.
It follows that $\zeta(t)< 2CE_0$ for all $t\ge 0$, yielding
\eqref{eq:cest} by definitions \eqref{czeta} and \eqref{Edef}.
\end{proof}

\br\label{diffusivermk}
\textup{
The rate of decay \eqref{eq:cest} is that of a heat kernel-- that
is, the mechanism for stability is {\it diffusive} only.
}
\er

\section{Stability of viscous shock solutions}\label{s:shock}

We turn now to the stability of viscous shock solutions of 
\eqref{burgers}, without
loss of generality, restricting to the case 
$$
\bar u(x)=-\tanh(x/2)
$$
described in \eqref{tanhform}.
Letting $\tilde u$ as before be a second solution of \eqref{burgers}, 
define the perturbation 
\be\label{spert}
u(x,t):=\tilde u(x+\alpha(t),t)-\bar u(x)
\ee
as the difference between a {\it translate of} $\tilde u$ and
the background wave $\bar u$, where the translation $\alpha(t)$
is to be determined later.

This yields after a brief computation the {\it perturbation equation}
\be\label{sperteq}
u_t-Lu= N(u)_x + \dot \alpha(t)(\bar u_x + u_x),
\ee
where 
$Lu:= u_{xx} -(a(x)u)_x $
is the linearization of $u_{xx}-f(u)_x$ about solution 
$\bar u=-\tanh(x/2)$, $a(x):=df(\bar u)(x)=\bar u(x)$,
and $N(u):= -u^2/2$ is the same quadratic order remainder
as in the constant-coefficient case.

\subsection{Linear solution operator/decomposition of 
the Green function}

The homogeneous linearized equation
\be\label{slin}
v_t-Lv=v_t+(a(x)v)_x-v_{xx}=0,
\quad v|_{t=0}=f
\ee
can be solved explicitly by linearized 
Hopf--Cole transformation \cite{S,N,Z3,GSZ}, 
to give an exact solution formula
\be\label{slinsol}
e^{Lt}f=\int_{-\infty}^{+\infty} G(x,t;y)f(y)dy,
\ee
where 
\ba
G(x,t;y):=e^{Lt}\delta_y(x) &=
\bar u'(x) 
\Big(\frac{1}{2}\Big)
\Big( \errfn (\frac{x-y-t}{\sqrt{4t}}) - \errfn (\frac{x-y+t}{\sqrt{4t}}) \Big) \\
&\quad
+
\Big( \Big(\frac{e^{-\frac{x}{2}}}{e^{\frac{x}{2}} + e^{-\frac{x}{2}}}\Big) 
\frac {e^{-\frac{(x-y-t)^2}{4t}} }{\sqrt{4 \pi t}} 
+ \Big(\frac{e^{\frac{x}{2}}} {e^{\frac{x}{2}} + e^{-\frac{x}{2}}} \Big) 
\frac { e^{-\frac{(x-y+t)^2} {4t}} } {\sqrt{4 \pi t}} \Big) \\
\ea
is the Green function for \eqref{slin} and
${\errfn} (z) := \frac{1}{2\pi} \int_{-\infty}^z e^{-\xi^2} d\xi$.

Following the approach of \cite{Z1,MaZ2,MaZ4}, decompose now
\be\label{Gdecomp}
G(x,t;y):=E(x,t;y)+ S(x,t;y)+ R(x,t;y),
\ee
where
\ba\label{E}
E(x,t;y)&:=\bar u'(x)e(y,t),
\qquad
e(y,t):= \Big(\frac{1}{2}\Big) 
\Big( \errfn (\frac{-y-t}{\sqrt{4t}}) - \errfn (\frac{-y+t}{\sqrt{4t}}) \Big), 
\ea
\be\label{S}
S(x,t;y):=
\Big( \Big(\frac{e^{-\frac{x}{2}}}{e^{\frac{x}{2}} + e^{-\frac{x}{2}}}\Big) 
\frac {e^{-\frac{(x-y-t)^2}{4t}} }{\sqrt{4 \pi t}} 
+ \Big(\frac{e^{\frac{x}{2}}} {e^{\frac{x}{2}} + e^{-\frac{x}{2}}} \Big) 
\frac { e^{-\frac{(x-y+t)^2} {4t}} } {\sqrt{4 \pi t}} \Big), 
\ee
and
\ba\label{R}
R(x,t;y)&:=
\bar u'(x) \Big(\frac{1}{2}\Big)
\Big( \errfn (\frac{x-y-t}{\sqrt{4t}}) - 
\errfn (\frac{-y-t}{\sqrt{4t}}) \Big) \\
&-\quad
\bar u'(x) \Big(\frac{1}{2}\Big)
\Big( \errfn (\frac{x-y+t}{\sqrt{4t}}) - 
\errfn (\frac{-y+t}{\sqrt{4t}}) \Big). \\
\ea
Here, the ``excited term'' $E$ represents the nondecaying part of
the linearized solution $v$, involving the zero-eigefunction
$L\bar u'=0$ associated with instantaneous translation of the 
background wave, the ``scattering term'' $S$ comprises Gaussian
signals convected along hyperbolic characteristics, and the
``remainder term'' $R$ a faster-decaying residual.

A straightforward calculation gives
\ba\label{Rcalc}
|R(x,t;y| &\le
C|x||\bar u'(x)|\int_0^1
\Big(\frac {e^{-\frac{(\theta x-y-t)^2}{4t}} }{\sqrt{4 \pi t}} 
-\frac {e^{-\frac{(\theta x-y+t)^2}{4t}} }{\sqrt{4 \pi t}} \Big)
d\theta\\
&\le
Ce^{-\theta |x|} \int_0^1
\Big(\frac {e^{-\frac{(\theta x-y-t)^2}{4t}} }{\sqrt{4 \pi t}} 
+\frac {e^{-\frac{(\theta x-y+t)^2}{4t}} }{\sqrt{4 \pi t}} \Big)
d\theta,
\ea
$\theta>0$, showing that $R$, as the product of an exponentially
decaying term and the sum of convected Gaussians, is indeed faster-decaying
than either $E$ or $S$.

\bl\label{Rlem1}
For some $C>0$, $\theta>0$, all $t>0$,
\ba\label{pR}
|R(x,t;y|&\le
Ce^{-\theta |x|/C}
\Big(\frac {e^{-\frac{( x-y-t)^2}{4Ct}} }{\sqrt{ t}} 
+\frac {e^{-\frac{(x-y+t)^2}{4Ct}} }{\sqrt{ t}} \Big)
+
Ce^{-\theta(|x-y|+t)},\\
|R_y(x,t;y|&\le
Ce^{-\theta |x|/C}
\Big(\frac {e^{-\frac{( x-y-t)^2}{4Ct}} }{ t} 
+\frac {e^{-\frac{(x-y+t)^2}{4Ct}} }{ t} \Big)
+
\frac{Ce^{-\theta(|x-y|+t)}}{\sqrt{t}}.\\
\ea
\el

\begin{proof}
Applying the Cauchy--Schwarz inequality in the argument of the
exponential, we find readily that
\be\label{temp}
e^{-\theta |x|/2} 
\Big(e^{-\frac{( \theta x-y-t)^2}{4t}}
+e^{-\frac{(\theta x-y+t)^2}{4t}} \Big)
\le
e^{-\theta |x|/C} \Big(e^{-\frac{( \theta x-y-t)^2}{4Ct}} 
+e^{-\frac{(x-y+t)^2}{4Ct}} \Big)
\ee
for $|x|\le Mt$ and $C>0$ sufficiently large, 
hence \eqref{Rcalc} implies \eqref{pR}(i). 
For $|y|>> |x|+|t$, \eqref{temp} holds trivially, likewise
giving \eqref{pR}(i).
In both of these cases, the lefthand side is bounded by the
first, Gaussian, term alone on the righthand side. 
In the remaining case $|x|>> t$ and $|y|\le M|x|$, 
we have for $C>0$ sufficiently large that
 $e^{-\theta|x|}\le e^{-(\theta/C)(|x-y|+t)}$,
from which we find directly from \eqref{R} that the lefthand side of 
\eqref{pR}(i) is bounded by the final term on the righthand side.

Similar computations yield \eqref{pR}(ii).
\end{proof}

\br\label{consrmk}
\textup{
The excited term $E$ converges as $t\to +\infty$
to $\bar u'(x)$ times 
$$
-\sigma(+\infty) := \int_{-\infty}^{+\infty}e(y,+\infty)f(y)dy=
(1/2)\int_{-\infty}^{+\infty} f(y)dy,
$$
the time-asymptotic state of the linearized equations \eqref{slin}
determined by conservation of mass (equals total integral
$\int_{-\infty}^{+\infty}v(x,t)dx$).
Note that $\bar u'(x)$ corresponds to infinitesimal translation
of the background wave $\bar u(x)$, hence a linear time-asymptotic
state $-\sigma \bar u'(x)$ corresponds roughly to a steady-state
perturbation $\bar u(x-\sigma)-\bar u(x)$ consisting of a shift,
or translation, $\sigma$ of the background wave.
The term $\sigma(t):=-\int_{-\infty}^{+\infty} e(y,t)f(y)dy$ 
thus measures, at a linearized level, the shift in location
of the shock at time $t$, or ``instantaneous shock shift''.
This refines the picture of behavior given by the time-asymptotic
shock shift $\sigma(+\infty)$.
%
}
\er

\bpr\label{greenbds}
The Green function $G$ decomposes as
$G=E+\tilde G$, $E=\bar u'(x)e(y,t)$,
where, for some $C>0$, all $t>0$,
\ba\label{sheatbds}
\Big|\int_{-\infty}^{+\infty} \tilde G(x,t;y)f(y)dy\Big|_{L^p(x)}&\le 
Ct^{-\frac{1}{2}(1-1/p)}|f|_{L^1},\\
\Big|\int_{-\infty}^{+\infty} \tilde G_y(x,t;y)f(y)dy\Big|_{L^p(x)}&\le 
Ct^{-\frac{1}{2}(1-1/p)-\frac{1}{2}}|f|_{L^1}.
\\
\ea
\ba\label{sheatbds2}
\Big|\int_{-\infty}^{+\infty} \tilde G(x,t;y)f(y)dy\Big|_{L^p(x)}&\le 
C |f|_{L^p},
\quad
\Big|\int_{-\infty}^{+\infty} \tilde G_y(x,t;y)f(y)dy\Big|_{L^p(x)}\le 
Ct^{-\frac{1}{2}}|f|_{L^p};\\
\ea
\ba\label{etbds}
\Big|\int_{-\infty}^{+\infty} e_t(y,t)f(y)dy\Big|&\le
Ct^{-\frac{1}{2}}|f|_{L^1},\qquad
\Big|\int_{-\infty}^{+\infty} e_{ty}(y,t)f(y)dy\Big|\le
Ct^{-1}|f|_{L^1},
\ea
\ba\label{eytbds}
\Big|\int_{-\infty}^{+\infty} e_{t}(y,t)f(y)dy\Big|&\le 
C |f|_{L^\infty},
\quad
\Big|\int_{-\infty}^{+\infty} e_{yt}(y,t)f(y)dy\Big|\le 
Ct^{-\frac{1}{2}}|f|_{L^\infty};\\
\ea
and
\ba\label{ebds}
\Big|\int_{-\infty}^{+\infty} e(y,t)f(y)dy\Big|&\le
C |f|_{L^1},\qquad
\Big|\int_{-\infty}^{+\infty} e_y(y,t)f(y)dy\Big|\le
Ct^{-1/2}|f|_{L^1}.
\ea
\epr

\begin{proof}
Defining $\tilde G:=R+S$, we have the decomposition $G=E+\tilde G$.
By \eqref{S} and estimate \eqref{pR}, $\tilde G$ and $\tilde G_y$
obey essentially the same bounds as $G$ and $G_y$ in the constant-coefficient
case \eqref{lGreen}, up to a harmless exponential error (the final
terms on the righthand sides of \eqref{pR}).
Thus, bounds \eqref{sheatbds} and \eqref{sheatbds2}
follow by the same argument used to prove \eqref{heatbds} and \eqref{heatbds2}.
By \eqref{E}, $|e_t|$ and $|e_{yt}|$ satisfy essentially the same bounds
as $\sup_x |\tilde G|$ and $\sup_x |\tilde G_y|$, hence
\eqref{etbds} and \eqref{eytbds} follow again from this same argument
in case $p=\infty$, which amounts to H\"older's inequality together
with $L^p$ bounds on $e$ and derivatives
(see Lemma \ref{3}, Appendix \ref{s:e} for a careful derivation
of these $L^p$ bounds).
Finally, \eqref{ebds} follows by $|e|\le C$, $|e_y|\le Ct^{-1/2}$
using the triangle inequality.
\end{proof}

\br\label{stabrmk}
\textup{
The apparently special Proposition \ref{greenbds} in fact
holds for viscous shock waves
of {\it general strictly parabolic systems}
provided that the shock satisfies
a generalized spectral stability, i.e., Evans function, condition
\cite{Z1,Z4,MaZ3}.
Indeed, there is a parallel decomposition of the Green function
as the sum of terms $E$, $S$, and $R$ 
with pointwise descriptions generalizing 
those 
of \eqref{E}, \eqref{S}, \eqref{pR}.
Similar bounds hold for Evans stable shocks of
general hyperbolic--parabolic systems \cite{MaZ3,Z4}.
Scalar shock waves are always spectrally stable, by the maximum
principle; hence, the stability condition does not make itself
apparent for Burgers equation.
}

\textup{
In the derivation of bounds by inverse Laplace transform
estimates, the terms $E$ and $S$ arise in a very natural
way as leading terms of a low-frequency ``scattering'' expansion
\cite{MaZ3,Z2} of the resolvent kernel about frequency $\lambda=0$,
{\it without the need to re-arrange terms as done here in 
the Burgers case}.
See Section 2, \cite{BeSZ}, for a particulary clear discussion
of the method from more general point of view.
Indeed, the decomposition of $G$ into $E$ and $\tilde G$
was suggested from the inverse Laplace transform point of
view \cite{ZH,Z1,MaZ3}.  Here, for pedagogical purposes,
we have imposed this structure by force on the explicit Green
function given by Hopf--Cole transformation in order to
demonstrate clearly the approach. 
}
\er

\subsection{Integral representation/$\alpha$-evolution scheme}

Recalling that $\bar u'(x)$ is a stationary
solution of the linearized equations $u_t=Lu$,
so that $L\bar u_x =0$, or
$$
\int^\infty_{-\infty}G(x,t;y)\bar u_x(y)dy=e^{Lt}\bar u_x(x)
=\bar u_x(x),
$$
we have, applying Duhamel's principle to \eqref{sperteq},
\be\label{prelim}
\begin{array}{l}
  \displaystyle{
  u(x,t)=\int^\infty_{-\infty}G(x,t;y)u_0(y)\,dy } \\
  \displaystyle{\qquad
  -\int^t_0 \int^\infty_{-\infty} G_y(x,t-s;y)
  (N(u)+\dot \alpha u ) (y,s)\,dy\,ds + \alpha (t) \bar u'(x).}
\end{array}
\ee
Defining $\alpha$ implicitly as
\begin{equation}
 \begin{array}{l}
  \displaystyle{
  \alpha (t)=-\int^\infty_{-\infty}e(y,t) u_0(y)\,dy }\\
  \displaystyle{\qquad
  +\int^t_0\int^{+\infty}_{-\infty} e_{y}(y,t-s)(N(u)+
  \dot \alpha\, u)(y,s) dy ds, }
  \end{array}
 \label{alpha}
\end{equation}
following \cite{ZH,Z4,MaZ2,MaZ3}, 
where $e$ is defined as in \eqref{E}, 
and substituting in \eqref{prelim} the decomposition $G=\bar u'(x)e +  \tilde G$ of Proposition 
\ref{greenbds}, 
we obtain the {\it integral representation}
\begin{equation}
\begin{array}{l}
 \displaystyle{
  u(x,t)=\int^\infty_{-\infty} \tilde G(x,t;y)u_0(y)\,dy }\\
 \displaystyle{\qquad
  -\int^t_0\int^\infty_{-\infty}\tilde G_y(x,t-s;y)(N(u)+
  \dot \alpha u)(y,s) dy \, ds, }
\end{array}
\label{u}
\end{equation}
and, differentiating (\ref{alpha}) with respect to $t$,
and observing that 
$e_y (y,s)\to 0$ as $s \to 0$, as the difference of 
approaching heat kernels,
\begin{equation}
 \begin{array}{l}
 \displaystyle{
  \dot \alpha (t)=-\int^\infty_{-\infty}e_t(y,t) u_0(y)\,dy }\\
 \displaystyle{\qquad
  +\int^t_0\int^{+\infty}_{-\infty} e_{yt}(y,t-s)(N(u)+
  \dot \alpha u)(y,s)\,dy\,ds. }
 \end{array}
\label{alphadot}
\end{equation}

Equations \eqref{u}, \eqref{alphadot} together form a complete
system in the variables $(u,\dot \alpha)$, 
from the solution of which we may afterward recover the
shift $\alpha$ via \eqref{alpha}.
From the original differential equation \eqref{sperteq}
together with \eqref{alphadot},
we readily obtain short-time existence and continuity with
respect to $t$ of solutions
$(u,\dot \alpha)\in L^1\cap L^\infty\times \RR$ 
by a 
standard
 contraction-mapping argument.\footnote{
Specifically, for initial time $T\ge 0$, and $t\ge T$, 
split the expression \eqref{alphadot} for $\dot \alpha(t)$ into the sum of
a bounded ``known'' term
  $-\int^\infty_{-\infty}e_t(y,t) u_0(y)\,dy 
  +\int^T_0\int^{+\infty}_{-\infty} e_{yt}(y,t-s)(N(u)+
  \dot \alpha u)(y,s)\,dy\,ds
$
and an 
``unknown term''
$ \int^t_T\int^{+\infty}_{-\infty} e_{yt}(y,t-s)(N(u)+ 
\dot \alpha u)(y,s)\,dy\,ds $ that is contractive 
for $(u,\dot \alpha)$ bounded and $|t-T|<<1$.
The $u$-equation \eqref{sperteq} may be treated in 
standard fashion, treating the righthand side as a forcing term and
expressing $u$ as an integral on $[T,t]$, again contractive
for 
$|t-T|<<1$.
}

\br\label{optrack}
\textup{
Here, the key step in deriving \eqref{u} is to observe that
the contribution in the righthand side of \eqref{prelim} coming
from terms involving $\bar u'(x)e(y,t)$ is, 
under definition \eqref{alpha}, exactly
$-\bar u'(x)\alpha(t)$, so cancels the final term.
That is, we have defined the instantaneous translation $\alpha(t)$
{\it from considerations of technical convenience}
so as to cancel all nondecaying terms in \eqref{prelim}.
Note that $\alpha(t)$ agrees to linear order with
the prescription $\sigma(t)$ in Remark \ref{consrmk} of
the instantaneous shock shift for the linearized equations.
}
\er

\subsection{Nonlinear iteration}

Associated with the solution $(u,\dot \alpha)$ of integral system 
\eqref{u}--\eqref{alphadot}, define 
\be\label{szeta}
\zeta(t):=\sup_{0\le s\le t,\, 1\le p\le \infty} 
\big( |u|_{L^p}(s)(1+t)^{\frac{1}{2}(1-1/p)}
+ |\dot \alpha(s)|(1+s)^{1/2} \big).
\ee

\bl\label{sclaim}
For all $t\ge 0$ for which $\zeta(t)$ is finite, some $C>0$,
and $E_0:=|u_0|_{L^1\cap L^\infty}$,
\be\label{eq:sclaim}
\zeta(t)\le C(E_0+\zeta(t)^2).
\ee
\el

\begin{proof}
With the established bounds on $\tilde G$ and $e$, the proof
of \eqref{eq:sclaim}
is almost identical to that of \eqref{eq:cclaim} in
the constant-coefficient case.
Noting, by quadratic dependence $N(u)=O(|u|^2)$ and the definition
\eqref{czeta} of $\zeta$, that 
\ba\label{sNbds}
|N(u)+\dot\alpha u|_{L^1}&\le 
C|u|_{L^1}(|u|_{L^\infty}+|\alpha|)\le \zeta(t)^2(1+t)^{-\frac{1}{2}}\\
|N(u)+\dot\alpha u|_{L^p}&\le C|u|_{L^p}(|u|_{L^\infty}+|\alpha|) \le 
\zeta(t)^2(1+t)^{-\frac{1}{2}(1-1/p)-\frac{1}{2}} ,
\ea
we obtain, similarly as in \eqref{cest},
applying Lemmas \ref{sheatbds}--\ref{sheatbds2} to representation 
\eqref{u}, 
\ba\label{sest}
|u(\cdot,t)|_{L^p(x)}& \le
C(1+t)^{-\frac{1}{2}(1-1/p)}E_0 +
C\zeta(t)^2\int_0^{t/2} (t-s)^{-\frac{1}{2}(1-1/p)-1/2}
(1+s)^{-\frac{1}{2}}ds\\
&\quad
+C\zeta(t)^2\int_{t/2}^{t} (t-s)^{-\frac{1}{2}}
(1+s)^{-\frac{1}{2}(1-1/p)-\frac{1}{2}}ds\\
&
\le
 C(E_0+\zeta(t)^2) (1+t)^{-\frac{1}{2}(1-1/p)}.
\ea
Similarly, by \eqref{etbds} and \eqref{eytbds},
\ba\label{sestad}
|\dot \alpha(t)|& \le
C(1+t)^{-\frac{1}{2}}E_0 +
C\zeta(t)^2\int_0^{t/2} (t-s)^{-1}
(1+s)^{-\frac{1}{2}}ds
\\ &\quad
+C\zeta(t)^2\int_{t/2}^{t} (t-s)^{-\frac{1}{2}}
(1+s)^{-1}ds\\
&
\le
C(E_0+\zeta(t)^2) (1+t)^{-\frac{1}{2}}.
\ea
Combining and rearranging \eqref{sest}--\eqref{sesta}, 
we obtain \eqref{cclaim}.
\end{proof}

\bc[Stability of shock solutions]\label{scor}
Viscous shock solutions $\bar u(x)$ of \eqref{burgers} are nonlinearly 
stable in $L^1\cap L^\infty$ and nonlinearly orbitally asympotically
stable in $L^p$, $p>1 $, with respect to 
initial perturbations $u_0$ that are sufficiently small in $L^1\cap L^\infty$.
More precisely, for some $C>0$ and $\alpha \in W^{1,\infty}(t)$, 
\ba\label{eq:smallsest}
|\tilde u-\bar u(\cdot -\alpha)|_{L^p}(t)&\le C(1+t)^{-\frac{1}{2}(1-1/p)}
|\tilde u-\bar u|_{L^1\cap L^\infty}|_{t=0},\\
|\dot\alpha(t)|&\le C(1+t)^{-\frac{1}{2}}
|\tilde u-\bar u|_{L^1\cap L^\infty}|_{t=0},\\
|\alpha(t)|&\le C
|\tilde u-\bar u|_{L^1\cap L^\infty}|_{t=0},\\
|\tilde u-\bar u|_{L^1\cap L^\infty}(t)&\le C 
|\tilde u-\bar u|_{L^1\cap L^\infty}|_{t=0},\\
\ea
for all $t\ge 0$, $1\le p\le \infty$, 
for solutions $\tilde u$ of \eqref{burgers} with
$|\tilde u-\bar u|_{L^1\cap L^\infty}|_{t=0}$ sufficiently small.
\ec

\begin{proof}
The first two inequalities follow by a proof identical to that
of Proposition \ref{ccor} in
the constant-coefficient case, using \eqref{sclaim}
and continuity of $\zeta$ wherever $\zeta$ is finite,
a consequence of short-time existence theory, to obtain
$\zeta(t)\le 2CE_0$, for
$E_0:=|\tilde u-\bar u|_{L^1\cap L^\infty}|_{t=0}\le \eta_0$ 
sufficiently small.
This yields the first two bounds by definition of $\zeta$.
The third then follows 
using \eqref{ebds}, by
\ba\label{sesta}
|\alpha(t)|& \le
C E_0 +
C\zeta(t)^2\int_0^{t/2} (t-s)^{-\frac{1}{2}}
(1+s)^{-\frac{1}{2}}ds
\\ &\quad
+C\zeta(t)^2\int_{t/2}^{t} (t-s)^{-\frac{1}{2}}
(1+s)^{-\frac{1}{2}}ds\\
&
\le
C(E_0+\zeta(t)^2). 
\ea

Finally, we note that 
$$
\tilde u(x,t)-\bar u(x)= 
u(x-\alpha(t),t)+ (\bar u(x)-\bar u(x-\alpha(t)),
$$
so that $|\tilde u(\cdot, t)-\bar u|$ is controlled
by the sum of $|u|$ and 
$|\bar u(x)-\bar u(x-\alpha(t))|\sim \alpha(t)|\bar u'(x)|$,
hence, by our estimates, remains $\le CE_0$ for all $t\ge 0$,
for $E_0$ sufficiently small.
This verifies the fourth inequality, yielding
 nonlinear stability and completing the result.
\end{proof}

\br\label{itmk}
\textup{
In the semilinear case considered here,
Corollary \ref{scor} could be proved in more straightforward
fashion by a contraction mapping argument applied directly to
the system \eqref{u}--\eqref{alphadot}, bypassing
the continuous induction argument above.
However, in more delicate situations such as  
the quasilinear parabolic or hyperbolic--parabolic
case, it is advantageous for reasons of regularity
to separate the issues of short-time existence/well-posedness
and long-time bounds, as we have done here; see \cite{MaZ2,MaZ4,Z4,RZ}
for further discussion.
}
\er

\br\label{sdiffusivermk}
\textup{
Again, the rate of decay \eqref{eq:sest} is that of a heat kernel-- that
is, the mechanism for stability is {\it diffusive} only, and not
involving compressivity of the shock.
This rate is in fact sharp, as may be seen intuitively
by considering a compactly supported perturbation
supported arbitrarily far from the shock location $x=0$.
Far from the shock, the background solution $\bar u$ is
approximately constant, and so behavior is like that of a perturbation
of a constant solution as studied in Section \ref{s:const}.
But, this is readily seen to decay like a heat kernel,
giving the stated rate \eqref{eq:sest}.
}
\er

\subsection{Postscript: phase-asymptotic vs. asymptotic orbital stability}
\label{s:phase}

A stronger condition that nonlinear orbital stability, proved above,
is {\it nonlinear phase-asymptotic orbital stability}, in which a 
perturbed solution $\tilde u$ is required to approach
not only the set of translates of $\bar u$,
but a {\it specific} translate of $u$.
In the language of Corollary \ref{scor}, this amounts to the requirement
that $\alpha$ have a limit
$\alpha(t)\to \alpha(+\infty)$ as $t\to +\infty$.

We do not establish this property in Corollary \ref{scor},
nor is it established in \cite{Z1,MaZ2,MaZ4}.
Indeed, for the general class of perturbations considered here
(and in \cite{Z1,MaZ2,MaZ4}), $\alpha(t)$ {\it if it converges to a limit}
does not do so at any uniform algebraic rate depending only on
$E_0$, $t$, as may be seen by considering perturbations with
support arbitrarily far from the shock location $x=0$.
See \cite{Z1} for further discussion. 

It is a strength of this approach that such data may be 
treated nonetheless, and in a simple fashion parallel to the
treatment of the constant-coefficient case.
However, phase-asymptotic stability does not seem to be
accessible by this simple argument scheme.
For proofs of phase-asymptotic stability under strengthened assumptions
on the initial data, involving additional pointwise information
on the solution, see \cite{R,HZ,HR,HRZ,RZ}.



\section{The system case}\label{s:sys}
We have described the nonlinear stability argument of \cite{Z1,MaZ2,MaZ4}
in the simple scalar setting of Burgers equation.
We now discuss briefly how this carries over to the case
of general hyperbolic-parabolic systems, including Navier--Stokes
equations of compressible gas dynamics and MHD.
Namely, Remark \ref{stabrmk} plus {\it essentially the same argument
described here} gives 
nonlinear orbital stability of viscous shocks provided that
they satisfy an {\it Evans function } (generalized spectral
stability) assumption yielding the necessary pointwise bounds. 
The Evans condition is necessary for linearized stability
as shown in \cite{ZH,MaZ3}.
It holds always for small-amplitude shocks, but may fail
in general for large-amplitude shocks.
In the large-amplitude case, it is readily checked numerically;
in certain special limits,
it may be checked analytically using asymptotic ODE
and or singular perturbation theory.
When the Evans condition fails, there are interesting implications
for dynamics/bifurcation; see \cite{Z5,Z7,TZ1,TZ2,TZ3,TZ4,SS,BeSZ}.

See \cite{AGJ,GZ,ZH,ZS,MaZ3} for discussion of the Evans function
and its origins.
For verification of the Evans condition for small-amplitude shocks,
see \cite{ZH,HuZ,PZ,FS1}.
For examples of unstable shocks, see \cite{GZ,ZS}.
For numerical and analytical 
verification for large-amplitude shocks, see \cite{BHZ,BHRZ,HLZ,HLyZ,CHNZ};
see \cite{Br,BrZ,BDG,HuZ2} for more general discussion of
numerical Evans function techniques.
See \cite{ZS,Z2,Z4,Z6,GMWZ,GMWZ2,FS2} for extensions to multiple dimensions.

The derivation of pointwise Green function bounds for general
systems is complicated,
involving detailed estimates on the resolvent kernel using Evans function
and asymptotic ODE
techniques, converted to bounds on the Green kernel
via stationary phase estimates in the inverse Laplace
transform formula.
See \cite{ZH,Z3,Z2,Z4,BeSZ,GMWZ,GMWZ2} for discussions of these 
and related techniques.
These are details of the {\it linear theory}.  Here, we have
chosen to isolate the {\it nonlinear
iteration argument} by restricting to a case (Burgers equation)
for which the linear theory is explicitly known a priori, in 
order to give the reader a flavor of the arguments.

We emphasize: {\it once the linearized theory is established, the nonlinear
shock-tracking argument of \cite{Z1,MaZ2,MaZ4,Z2} 
is essentially the same for system or for scalar case}.
See Remark \ref{stabrmk}.

\medskip
{\bf Acknowledgement.}  Thanks to Mark Williams and Benjamin
Texier for their interest in the work,
and for several helpful comments improving the exposition.

\bigskip
{\bf APPENDICES}

\appendix
\section{The small-amplitude limit}
It is instructive to consider the small-amplitude limit
$|u_+-u_-|\to 0$.
Consider now the {family} of stationary viscous shock solutions
\be\label{fam}
\bar u^\eps (x):= -\eps \tanh(\eps x/2),
\qquad
\lim_{x\to \pm \infty}= \bar u^\eps(x) u_\pm^\eps=\mp \eps
\ee
of \eqref{burgers}, and examine behavior as $\eps \to 0$.

Denote the associated homogeneous linearized equation by
\be\label{smallslin}
v_t-L^\eps v=v_t+(a^\eps (x)v)_x-v_{xx}=0,
\quad v|_{t=0}=f
\ee
where $a^\eps(x):=\bar u^\eps(x)$.
The invariance $(x,t,u)\to (x/\eps , t/\eps^2 , u/\eps)$ of
Burgers equation converts this to the $\eps$-independent case \eqref{tanhform}
considered in Section \ref{s:shock}, from which we may deduce
the $\eps$-dependent Green function formula
\be\label{smallslinsol}
e^{L^\eps t}f=\int_{-\infty}^{+\infty} G^\eps(x,t;y)f(y)dy,
\ee
where 
\ba
G^\eps(x,t;y):=e^{L^\eps t}\delta_y(x) &=
(\bar u^\eps)'(x) 
\Big(\frac{1}{2\eps }\Big)
\Big( \errfn (\frac{x-y-\eps t}{\sqrt{4t}}) - \errfn (\frac{x-y+\eps t}{\sqrt{4t}}) \Big) \\
&\quad
+
\Big( \Big(\frac{e^{-\frac{\eps x}{2}}}{e^{\frac{\eps x}{2}} + e^{-\frac{\eps x}{2}}}\Big) 
\frac {e^{-\frac{(x-y-\eps t)^2}{4t}} }{\sqrt{4 \pi t}} 
+ \Big(\frac{e^{\frac{\eps x}{2}}} {e^{\frac{\eps x}{2}} + e^{-\frac{\eps x}{2}}} \Big) 
\frac { e^{-\frac{(x-y+\eps t)^2} {4t}} } {\sqrt{4 \pi t}} \Big) \\
\ea
and $(\bar u^\eps)'(x)=\eps^2 \bar u'(\eps x)\sim
\eps^2 e^{-\theta \eps |x|}$, $\theta>0$.
Here, we are using the scaling relations
$$
\bar u^\eps(x)=\eps \bar u(\eps x)
\quad \hbox{\rm and }\quad
G^\eps(x,t;y)=\eps G(\eps x, \eps^2 t;\eps y).
$$

Decompose again
\be\label{smallGdecomp}
G^\eps(x,t;y):=E^\eps(x,t;y)+ S^\eps(x,t;y)+ R^\eps(x,t;y),
\ee
where
\ba\label{smallE}
E^\eps(x,t;y)&:=(\bar u^\eps)'(x)e^\eps(y,t),
\qquad
e^\eps(y,t):= \Big(\frac{1}{2\eps }\Big) 
\Big( \errfn (\frac{-y-\eps t}{\sqrt{4t}}) - \errfn (\frac{-y+\eps t}{\sqrt{4t}}) \Big), 
\ea
\be\label{smallS}
S^\eps (x,t;y):=
\Big( \Big(\frac{e^{-\frac{\eps x}{2}}}{e^{\frac{\eps x}{2}} + e^{-\frac{\eps x}{2}}}\Big) 
\frac {e^{-\frac{(x-y-\eps t)^2}{4t}} }{\sqrt{4 \pi t}} 
+ \Big(\frac{e^{\frac{\eps x}{2}}} {e^{\frac{\eps x}{2}} + e^{-\frac{\eps x}{2}}} \Big) 
\frac { e^{-\frac{(x-y+\eps t)^2} {4t}} } {\sqrt{4 \pi t}} \Big), 
\ee
and
\ba\label{smallR}
R^\eps (x,t;y)&:=
(\bar u^\eps) '(x) \Big(\frac{1}{2\eps }\Big)
\Big( \errfn (\frac{x-y-\eps t}{\sqrt{4t}}) - 
\errfn (\frac{-y-\eps t}{\sqrt{4t}}) \Big) \\
&+\quad
(\bar u^\eps)'(x) \Big(\frac{1}{2\eps}\Big)
\Big( \errfn (\frac{-y+\eps t}{\sqrt{4t}}) - 
\errfn (\frac{x-y+\eps t}{\sqrt{4t}}) \Big). \\
\ea

Defining the perturbation 
\be\label{smallpertu}
u(x,t):=\tilde u(x-\alpha(t),t)-\bar u^\eps(x),
\ee
setting $\tilde G^\eps:=S^\eps + R^\eps$, 
and following the steps of Section \ref{s:shock}, we obtain
again the { integral representation}
\begin{equation}
\begin{array}{l}
 \displaystyle{
  u(x,t)=\int^\infty_{-\infty} \tilde G(x,t;y)^\eps u_0(y)\,dy }\\
 \displaystyle{\qquad
  -\int^t_0\int^\infty_{-\infty}\tilde G^\eps _y(x,t-s;y)(N(u)+
  \dot \alpha u)(y,s) dy \, ds, }
\end{array}
\label{smallu}
\end{equation}
\begin{equation}
 \begin{array}{l}
 \displaystyle{
  \dot \alpha (t)=-\int^\infty_{-\infty}e^\eps_t(y,t) u_0(y)\,dy }\\
 \displaystyle{\qquad
  +\int^t_0\int^{+\infty}_{-\infty} e^\eps_{yt}(y,t-s)(N(u)+
  \dot \alpha u)(y,s)\,dy\,ds. }
 \end{array}
\label{smallalphadot}
\end{equation}
\begin{equation}
 \begin{array}{l}
  \displaystyle{
  \alpha (t)=-\int^\infty_{-\infty}e^\eps(y,t) u_0(y)\,dy }\\
  \displaystyle{\qquad
  +\int^t_0\int^{+\infty}_{-\infty} e^\eps_{y}(y,t-s)(N(u)+
  \dot \alpha\, u)(y,s) dy ds. }
  \end{array}
 \label{smallalpha}
\end{equation}

{\bf Dependence on $\eps$.}
Evidently, we could carry through the entire stability analysis
of Section \ref{s:shock}, as the $\eps$-dependent
Green function $G^\eps= E^\eps+S^\eps+ R^\eps$ has the same form as $G$.
However, the bounds obtained in this way-- in particular, the
estimate \eqref{pR} on the remainder $R$-- would involve constants
$C=C(\eps)>0$ blowing up as $\eps\to 0$.
This means that the allowable size $E_0\le \frac{1}{4C( \eps)^2}$
of perturbations, determined in the proof of
Corollary \ref{ccor}, goes to zero as $\eps\to 0$.
That is, {\it the basin of attraction of the shock $\bar u^\eps$ 
established by our basic stability argument shrinks to zero as $\eps\to 0$}.
Indeed, the bounds derived for general systems in \cite{ZH,MaZ3}
(described briefly in Section \ref{s:sys})
share this same property, and so the basin of attraction for
the stability results proved in \cite{MaZ2,MaZ4,Z1,Z2,HZ} and
related works go
to zero as the shock amplitude goes to zero.

%

However, this is not an inherent limitation of the method,
or the shock.
Following, we show that by different, more careful, estimates
of $E^\eps$ and $R^\eps$, we may in fact recover a uniform stability
result, valid for perturbations of sufficiently small size
{\it independent of $\eps$.}

\bpr\label{smallgreenbds}
For some $C>0$ {\bf independent of $\eps$}, $0<\eps \le 1$, and all $t>0$,
\ba\label{smallsheatbds}
\Big|\int_{-\infty}^{+\infty} \tilde G^\eps(x,t;y)f(y)dy\Big|_{L^p(x)}&\le 
Ct^{-\frac{1}{2}(1-1/p)}|f|_{L^1},\\
\Big|\int_{-\infty}^{+\infty} \tilde G^\eps_y(x,t;y)f(y)dy\Big|_{L^p(x)}&\le 
Ct^{-\frac{1}{2}(1-1/p)-\frac{1}{2}}|f|_{L^1}.
\\
\ea
\ba\label{smallsheatbds2}
\Big|\int_{-\infty}^{+\infty} \tilde G^\eps(x,t;y)f(y)dy\Big|_{L^p(x)}&\le 
C |f|_{L^p},
\quad
\Big|\int_{-\infty}^{+\infty} \tilde G^\eps_y(x,t;y)f(y)dy\Big|_{L^p(x)}\le 
Ct^{-\frac{1}{2}}|f|_{L^p};\\
\ea
\ba\label{smalletbds}
\Big|\int_{-\infty}^{+\infty} e^\eps_t(y,t)f(y)dy\Big|&\le
Ct^{-\frac{1}{2}}|f|_{L^1},\qquad
\Big|\int_{-\infty}^{+\infty} e^\eps_{ty}(y,t)f(y)dy\Big|\le
Ct^{-1}|f|_{L^1},
\ea
\ba\label{smalleytbds}
\Big|\int_{-\infty}^{+\infty} e^\eps_{t}(x,t;y)f(y)dy\Big|&\le 
C |f|_{L^\infty},
\quad
\Big|\int_{-\infty}^{+\infty} e^\eps_{yt}(x,t;y)f(y)dy\Big|\le 
Ct^{-\frac{1}{2}}|f|_{L^\infty};\\
\ea
and
\ba\label{smallebds}
\Big|\int_{-\infty}^{+\infty} e^\eps(y,t)f(y)dy\Big|&\le
C\eps^{-1} |f|_{L^1},\qquad
\Big|\int_{-\infty}^{+\infty} e^\eps_y(y,t)f(y)dy\Big|\le
C\eps^{-1}t^{-1/2}|f|_{L^1}.
\ea
\epr

\begin{proof}
As $S^\eps$ evidently obeys the same decay estimates as $S$,
to establish the stated bounds on $\tilde G^\eps$,
it is sufficient to establish them for $R^\eps$.
This is a straightforward consequence of Lemmas \ref{rt}
and \ref{Rbds} established in Appendix \ref{s:Re}.
Likewise, for the stated bounds on $e^\eps$ it is sufficient
to establish corresponding $L^p(y)$ bounds on $\eps^\eps$,
from which the results then follow by H\"older's inequality.
The needed bounds are established in Lemma \ref{sebds},
Appendix \ref{s:ee}.
\end{proof}

\bc[Stability of small-amplitude shock solutions]\label{smallscor}
For $0<\eps\le 1$,
viscous shock solutions $\bar u^\eps(x)$ of \eqref{burgers} are nonlinearly 
stable in $L^1\cap L^\infty$ and nonlinearly orbitally asympotically
stable in $L^p$, $p>1 $, with respect to 
initial perturbations $u_0$ with $L^1\cap L^\infty$ norm less than
or equal to $\eta_0>0$ sufficiently small, where $\eta_0$ is
{independent of $0<\eps\le 1$}.
More precisely, for some $C>0$ {independent of $0<\eps\le 1$}, 
there is $\alpha \in W^{1,\infty}(t)$ such that 
\ba\label{eq:sest}
|\tilde u-\bar u^\eps(\cdot -\alpha)|_{L^p}(t)&\le C(1+t)^{-\frac{1}{2}(1-1/p)}
E_0,\\
|\dot\alpha(t)|&\le C(1+t)^{-\frac{1}{2}}E_0,\\
|\alpha(t)|&\le C \eps^{-1}E_0,\\
|\tilde u-\bar u^\eps|_{L^1\cap L^\infty}(t)&\le C E_0,\\
\ea
for all $t\ge 0$, $1\le p\le \infty$, 
for solutions $\tilde u$ of \eqref{burgers} with
$E_0:=|\tilde u-\bar u^\eps|_{L^1\cap L^\infty}|_{t=0}\le \eta_0$.
\ec

\begin{proof}
The proof of the first two bounds follows exactly as in the proof 
of Corollary \ref{scor} in the fixed-amplitude case, 
since the integral equations for $(u,\dot \alpha)$ form a closed system
involving only  $\tilde G^\eps$, $e^\eps_t$ and $e^\eps_{yt}$,
and the bounds on $\tilde G^\eps$,
$e^\eps_t$ and $e^\eps_{yt}$ are the same as the bounds on
on $\tilde G$, $e_t$ and $e_{yt}$ in the fixed-amplitude case.
With these bounds established, we obtain the third bound from
\eqref{smallalpha}, using the fact that the bounds on $e^\eps$
and $e^\eps_y$ are no worse than $\eps^{-1}$ times the bounds
on $e$ and $e_y$ in the fixed-amplitude case.

Finally, we note that 
$
\tilde u(x,t)-\bar u^\eps(x)= 
u(x-\alpha(t),t)+ (\bar u^\eps(x)-\bar u^\eps(x-\alpha(t)),
$
so that $|\tilde u(\cdot, t)-\bar u^\eps|$ is controlled
by the sum of $|u|$ and 
$|\bar u^\eps(x)-\bar u^\eps(x-\alpha(t))|$.
By monotonicity of scalar shock profiles as orbits of
the first-order scalar profile ODE \eqref{prof},
$ \bar u^\eps(x)-\bar u^\eps(x-\alpha(t)) $ has one sign,
hence 
$$
 |\bar u^\eps(x)-\bar u^\eps(x-\alpha(t))|_{L^1}=
\Big|\int_{-\infty}^{+\infty}
(\bar u^\eps(x)-\bar u^\eps(x-\alpha(t)))dx \Big|
=
|\alpha(t)||u^\eps_+-u^\eps_-|,
$$
and, by \eqref{eq:sest}(iii),
$$
 |\bar u^\eps(x)-\bar u^\eps(x-\alpha(t))|_{L^1}=
2\eps |\alpha(t)| \le 2CE_0.
$$
Likewise, by the Mean Value Theorem,
$$
|\bar u^\eps(x)-\bar u^\eps(x-\alpha(t))|\le 
|\alpha(t)|(\bar u^\eps)'|_{L^\infty}
\le
(CE_0/\eps)(\eps^2)=CE_0 \eps,
$$
by the asymptotics $\bar u^\eps)'\sim \eps^2 e^{-\theta  \eps |x|}$.
Thus, 
$|\bar u^\eps(x)-\bar u^\eps(x-\alpha(t))|_{L^1\cap L^\infty}
\le CE_0$, and so
$|\tilde u(x,t)-\bar u^\eps(x)|_{L^1\cap L^\infty}\le CE_0$
for all $t\ge 0$, for $E_0$ sufficiently small.
This verifies the fourth inequality, yielding
 nonlinear stability and completing the result.
\end{proof}

\br\label{largeshift}
\textup{
In the small-amplitude limit $\eps\to 0$,
the shock shift $\alpha \to +\infty$
as $\eps^{-1}$ times perturbation mass.
Nonetheless, the stability estimates are uniform, independent of $\eps$.
}
\er

\br
\textup{
As discussed in Section \ref{s:phase},
we have obtained stability for a class $L^1\cap L^\infty$ 
of perturbations that lead to shock shifts $\alpha$ not only of
order $1/\eps$, but also decaying subalgebraically to their
limits $\alpha(+\infty)$, if they exist.
}
\er


\br\label{smallsys}
\textup{
Here we have treated only the simple and explicit case of
Burgers equation.
It would be very interesting to try to treat the small-amplitude
 system case by a similarly simple argument based on this
approach, using the singular perturbation techniques developed
in \cite{MaZ3,PZ} to obtain the necessary sharpened $\eps$-dependent
bounds analogous to those of Proposition \ref{smallgreenbds} 
in the Burgers case to try to obtain results uniform in $\eps$.
}
\er

\section{Alternative shock-tracking schemes}\label{s:alt}

As discussed in Remarks \ref{consrmk} and \ref{stabrmk},
the quantity $\alpha(t)$ introduced for technical
reasons in \eqref{alpha}, has an interpretation as an 
``instantaneous shock shift'', measuring the approximate location of
a perturbed viscous shock profile at time $t$.
This suggests the
question what is the ``exact'' location of an asymptotic shock profile,
and how well $\alpha(t)$ approximates this location.
The study of this question leads to an interesting
class of alternative shock-tracking
schemes that are time-asymptotically equivalent to 
\eqref{u}--\eqref{alphadot}, based on localized projections,
converging in the ``infinite-localization'' limit
to a pointwise phase condition introduced in \cite{GMWZ}
in the context of the small-viscosity limit.

Unlike a perturbed inviscid shock wave, which is sharply located
by the presence a discontinuity, a perturbed viscous shock wave is 
{\it smooth}, so requires some extrinsic criterion to define its location.
Two intuitive definitions immediately come to mind.  
The first, defining the location of an unperturbed stationary
scalar shock $u\equiv \bar u(x)$ without loss of generality
to be the origin, $x=0$,
is simply to define the location $\alpha(t)$ of a perturbed shock
$\tilde u$ as the point $\alpha(t)$ at which
$\tilde u$ takes on the value $\bar u(0)$, or
\be\label{scal}
\tilde u(\alpha(t),t)= \bar u(0).
\ee
By the Implicit Function Theorem and the fact that $\bar u'(0)\ne 0$
(recall that $\bar u$ is monotone, as the solution of a scalar first-order 
traveling-wave ODE), this uniquely defines $\alpha$ for
$|\tilde u'- \bar u'|_{L^\infty}(t)$ sufficiently small.

In the system case $u\in \RR^n$, we cannot satisfy \eqref{scal}
for all $n$ coordinates using the single parameter $\alpha$,
so we must choose some preferred coordinate direction,
substituting for \eqref{scal} the system analog
\be\label{1gmwz}
\ell \cdot \tilde u(\alpha(t),t)= \ell \cdot \bar u(0)
\ee
for some vector $\ell \in \RR^n$
such that $\ell \cdot \bar u'(0)\ne 0$,
a condition that, by the Implicit Function Theorem, guarantees that
$\alpha(t)$ is well-defined for 
$|\tilde u'- \bar u'|_{L^\infty}(t)$ sufficiently small.

Defining the perturbation variable
\be\label{pertrem}
u(x,t)=\tilde u(x+\alpha(t),t)-\bar u(x), 
\ee
following the notation of Section \ref{s:shock},
we find that
\eqref{1gmwz} translates to the {\it phase condition }
\be\label{gmwz}
 \ell \cdot u(0,t)=0,
\ee
determining $\alpha(t)$ implicitly through \eqref{pertrem}.
Condition \eqref{gmwz} is particularly natural from the point of view
of the resolvent equation arising in
solution by Laplace transform of the associated linearized equations.
For, the resolvent equation consists of an underdetermined 
ordinary differential boundary-value problem for which the 
standard treatment is to remove indeterminacy by one or more 
phase conditions like \eqref{gmwz}.
Indeed, this condition was introduced in \cite{GMWZ} starting from just
such considerations,
for the study of shock stability in the vanishing viscosity
limit,\footnote{
More precisely, a multi-dimensional version reducing to 
\eqref{gmwz} in the one-dimensional case.}

%

%
%
%
%


The second intuitive definition is, following Goodman \cite{G}, 
to define the shock shift $\alpha$ so as to minimize the least squares
distance of $\tilde u(x,t)$ from the shifted shock 
$\bar u(x-\alpha(t)$, that is, to minimize
$|u(\cdot,t)|_{L^2}$.
This leads to the ``localized projection condition''
(Euler-Lagrange equation)
\be\label{ls}
\langle \ell , u\rangle_{L^2}=0,
\qquad
\langle \ell, \bar u'\rangle_{L^2}=1,
\ee
where $\ell(x):=\frac{\bar u'(x)}{|\bar u'|_{L^2}^2}$
(see Appendix \ref{s:el} for this calculation).
Here, the word ``localized'' refers to the fact that
$\ell(x)$ decays as $x\to \pm \infty$.
More generally, we denote as a {localized projective condition}
any condition of form \eqref{ls} 
with $\ell\in L^1$.
This can be viewed as a nonlocal 
version of the pointwise phase condition \eqref{gmwz},
converging to \eqref{gmwz} in the ``infinite-localization limit''
$\ell(x)\to \ell_0 \delta(x)$, $\ell_0\in \RR^n$ constant,
of a Dirac measure.

Each of these schemes (either of form \eqref{gmwz} or \eqref{ls})
may be written as an evolution equation in $(u,\alpha)$.
Defining the perturbation variable $u$ of \eqref{pertrem}, we 
find as in Section \ref{s:shock} that $u$ obeys the partial differential
equation
\be\label{evrem}
u_t-Lu=N(u)_x +\dot \alpha (\bar u_x+u_x)
\ee
depending on $\dot \alpha$, defined implicitly by \eqref{ls}.
Differentiating \eqref{ls} with respect to $t$, we obtain
$$
0=\langle \ell , u_t \rangle_{L^2}=
\langle \ell, Lu +N(u)_x +\dot \alpha (\bar u_x+ u_x)\rangle_{L^2},
$$
which, using $\langle \ell, \bar u_x\rangle_{L^2}=1$, reduces to
$\dot \alpha(1+\langle \ell, u_x\rangle_{L^2}= 
-\langle \ell, Lu +N(u)_x \rangle_{L^2}$,
or, rearranging,
\be\label{impdef}
\dot \alpha=
- \frac{\langle \ell, Lu +N(u)_x \rangle_{L^2}}
{1+ \langle \ell, u_x\rangle_{L^2} },
\ee
well-defined for $u\in H^2$ with $|u|_{H^2}$ sufficiently small.
See \cite{G,TZ1,Z7} for related discussion.

Together, \eqref{evrem}--\eqref{impdef}
determine a closed system of evolution equations for
$(u,\dot\alpha)$, similar in spirit to the system
\eqref{u}--\eqref{alphadot} of Section \ref{s:shock},
but {\it local in time}, whereas 
the system \eqref{u}--\eqref{alphadot} 
involves ``memory terms''
depending on values of $u$, $\dot \alpha$ at earlier
times $ s\le t$.
For each choice of test function $\ell$, there results a different
evolution system, and different solutions $(u,\dot \alpha)$ and 
$\alpha$, 
representing different decompositions of the common solution $\tilde u$
of \eqref{burgers} under investigation, a perturbed viscous shock wave.

We know already from the analysis of  Corollary \ref{scor} that
the solution $\tilde u$ exists for all time, and converges to the
set of translates of the background shock $\bar u$.
However, it is not a priori clear that the system 
\eqref{evrem}--\eqref{impdef} has a global solution for any
particular choice of $\ell$, nor that the solution $u$ should
decay as $t\to 0$.
That is, it is not clear which of these alternative shock
tracking schemes 
gives an accurate estimate of shock location in the sense
that the known convergence of $\tilde u$ to the set of translates is 
revealed by decay at the appropriate rate of the perturbation
variable $u$.

The following proposition asserts that {\it all} of these schemes
are accurate in this sense, so that in principle any one of them
could be used as the basis of an argument for nonlinear stability.
Indeed, all lead to the same rates of decay.

\bpr\label{comparison}
Let $\uref$, $\alpharef$ denote the solution of \eqref{u}--\eqref{alphadot}
of Section \ref{s:shock}, with initial data
$\tilde u_0-\bar u$, $E_0:=
|\tilde u_0-\bar u |_{L^1\cap H^2}$ sufficiently small,
and $u$, $\alpha$ denote the solution with same initial data
of \eqref{evrem}--\eqref{impdef},
with $\ell\in L^1$.
Then, $u$, $\uref$ exist for all $t\ge 0$, with
\ba\label{decay}
|u|_{L^p}(t),\, |\uref|_{L^p}(t) &\le CE_0 (1+t)^{-\frac{1}{2}(1-1/p)},\\
|u|_{H^2}(t),\, |\uref|_{H^2}(t) &\le C E_0 (1+t)^{-1/4},\\
|u|_{L^1\cap H^2}(t)- |\uref|_{L^1\cap H^2}(t)&\le C E_0 (1+t)^{-1/2},\\
|\tilde u-\bar u|_{L^1\cap H^2}(t) &\le C E_0,\\
|\alpha|(t),\, |\alpharef|(t) &\le C E_0,\\
|\alpha- \alpharef|(t)&\le C E_0 (1+t)^{-1/2} .
\ea
\epr

\begin{proof}
A routine extension of the proof of Corollary \ref{scor},
using the additional assumption of $H^2$ smallness of the
initial data yields \eqref{eq:smallsest} augmented with
$ |\uref|_{H^2}(t)\le CE_0 (1+t)^{-\frac{1}{4}}, $
We omit the details. (But see the results of 
\cite{MaZ2,MaZ4} in the much more complicated
system case.)
The corresponding bounds \eqref{decay}(i)--(ii), hence global
existence of $u$, thus follow
provided that we can establish \eqref{decay}(iii).

Expanding
\ba\label{keycomp}
u(x,t)&=\tilde u(x+\alpha(t),t)-\bar u(x)\\
&=\tilde u(x+\alpha(t),t)
- \bar u(x+ (\alpha-\alpharef))
+ \bar u(x+ (\alpha-\alpharef)) - \bar u(x)  \\
&=
\uref(x+(\alpha-\alpharef),t)
+ \big(\bar u(x+ (\alpha-\alpharef)) - \bar u(x) \big) ,
\ea
we find using the Triangle inequality, followed
by the Mean Value Theorem together with exponential
decay of $\bar u'$, that
$$
\begin{aligned}
|u|_{L^1\cap H^2}(t)- |\uref|_{L^1\cap H^2}(t)&\le 
\big|\bar u(x+ (\alpha-\alpharef)) - \bar u(x) \big|_{L^1\cap H^2}
\le
C|\alpha-\alpharef|(t),
\end{aligned}
$$
so that \eqref{decay}(iii) follows from \eqref{decay}(vi).
Likewise, (iv) follows from (i)--(iii) and (v),
which in turn follows from (vi) and the bounds 
on $|\alpharef|(t)$ established
in Section \ref{s:shock}.

Thus, it remains only to prove \eqref{decay}(vi).
Applying definition $\langle \ell, u\rangle_{L^2}=0$ to
expansion \eqref{keycomp}, we obtain
\ba\label{compfpre}
\langle \ell, \uref(x+(\alpha-\alpharef),t)\rangle_{L^2}&=
-\langle \ell, \bar u(x+ (\alpha-\alpharef)) - \bar u(x) \rangle_{L^2}\\
&=
-\langle \ell,  (\alpha-\alpharef)) \bar u'+ O(
|\alpha-\alpharef)|^2 ) \rangle_{L^2}.\\
\ea
Applying now $\langle \ell, \bar u'\rangle_{L^2}=1$, and 
rearranging, we obtain
\ba\label{compf}
 |\alpha-\alpharef|(t)&\le
| \ell|_{L^1}
( |\uref|_{L^\infty}(t)+ C |\alpha-\alpharef|^2)\\
& \le C_2 ( E_0(1+t)^{-1/2} + |\alpha-\alpharef|^2),
\ea
yielding \eqref{decay}(vi) provided 
 $|\alpha-\alpharef|$ is sufficiently small.
The result then follows by continuity of $\alpha$, $\alpharef$
and smallness of $\alpha$ at $t=0$ for $E_0$ small, recalling
that $\alpharef(0)=0$.
\end{proof}

\br\label{measure}
\textup{
As the only bound used on $\ell$ was its $L^1$ norm,
the proof of Proposition \ref{comparison} is easily adapted to
the case that $\ell$ is a bounded measure,
in particular the case of a phase condition \eqref{gmwz}.
This includes also more general cases such as the sum of point measures,
leading to a sort of ``difference stencil'' condition determining
shock location.
}
\er

\br\label{trackrmk1}
\textup{
Recalling that $\alpharef(t)$ in general decays at most at subalgebraic
rate (see Remark \ref{sdiffusivermk}), we see from \eqref{decay}(iv) 
that $\alpha$ and $\alpharef$ are {\it time-asymptotically
equivalent} in the sense that
$|\alpha-\alpharef|$ decays at a rate faster than the
(general) rate of decay of $|\alpharef|$.
}
\er

\br\label{trackrmk2}
\textup{
For initial data in addition decaying as $|u_0(x)|\le CE_0(1+|x|)^{-3/2}$,
it is shown for general systems in \cite{HR,RZ}
that $\alpharef$ decays at the faster rate
\be\label{salphabd}
|\alpharef(t) -\alpharef(+\infty)|\le CE_0(1+t)^{-1/2}.
\ee
However, the same analysis yields sharpened bounds on $\uref$
as well, giving also
$$
|\uref(x,t)|\le CE_0(1+t)^{-1}
\; \hbox{\rm for }\; |x|\le \theta t,
$$
$\theta>0$ sufficiently small.
Substituting in \eqref{compfpre}, we obtain in place of 
\eqref{compf} the estimate
 $$
|\alpha-\alpharef|(t)\le
CE_0(1+t)^{-1} |\ell|_{L^1}
+CE_0 (1+t)^{-1/2}\int_{|x|\ge \theta t} |\ell(x)|dx
+C |\alpha-\alpharef|^2,
$$
yielding $|\alpha -\alpharef|(t)\le CE_0(1+t)^{-1}$ provided
$|\ell(x)|\le C(1+|x|)^{-3/2}$. 
}

\textup{
Thus, under this strengthened decay requirement on $\ell$,
we obtain time-asymptotic equivalence of $\alpha$ and $\alpharef$ 
also in this case.
Bound \eqref{salphabd} is sharp, as can be seen by direct computation
on the linear term in \eqref{alpha} for data decaying as $(1+|x|)^{-3/2}$.
(Note that the linear $O(E_0)$ term dominates the nonlinear $O(E_0^2)$
term up to any finite time, for $E_0$ sufficiently small.)
}
\er

{\bf Conclusions.}
By comparison with the scheme of Section \ref{s:shock}, we find
that each of the alternative shock-tracking schemes described
in this Appendix, based on localized phase conditions, yields
a globally defined solution exhibiting the same rates of decay
as the perturbation $\uref$ defined in Section \ref{s:shock}.
That is, essentially any tracking scheme based on information
that is ``local to the shock'' in the sense that it is accessible by inner
product with an $L^1$ function (resp. bounded measure) 
$\ell$ yields a convergent system of perturbation equations.
Note, further, that the only information used to draw these conclusions
consists of estimates on $(\uref,\alpharef)$ already established 
in \cite{Z1,MaZ2,MaZ3,MaZ4,HZ,RZ} for 
Evans-stable Lax or undercompressive type
shocks of general hyperbolic--parabolic systems. 
{\it Thus, the conclusions of Proposition \ref{comparison} and 
Remarks \ref{measure}--\ref{trackrmk2} remain
valid for Evans-stable Lax or undercompressive shocks of general systems of 
hyperbolic--parabolic conservation laws.}\footnote{
With the inclusion of additional
phase conditions to account for additional degrees of freedom
in the time-asymptotic state (see \cite{HZ,RZ}), 
these methods and estimates extend also 
to the overcompressive case.}

An interesting question is whether we could carry
out a nonlinear stability analysis for these schemes
{\it from first principles } rather than by comparison to
our existing results.
This is particularly intriguing for the case of the pointwise
phase condition \eqref{gmwz}, for which resolvent (and thus
pointwise Green function) bounds are available through
the framework developed in \cite{GMWZ}.
Besides the intrinsic interest of this question, there are
real advantages to the scheme based on \eqref{gmwz} for extension
to more complicated situations: for example, the fact that it
is local in time (the scheme in Section \ref{s:shock} by contrast
involves ``memory terms``), and that the phase condition \eqref{gmwz}
makes no reference to the explicit structure of the system.

\section{Miscellaneous estimates}

\subsection{Bounds on $e$}\label{s:e}

\bl\label{3} For some $C>0$ and all $t>0$,
\be\label{35}
|{e} (\cdot, t)|_{L^\infty}, \le C ,
\ee
\be\label{36}
|{e}_y (\cdot, t)|_{L^p},  |{e}_t(\cdot, t)|_{L^p} 
\le C t^{-\frac{1}{2}(1-1/p)},
\ee
\be\label{37}
|{e}_{ty}(\cdot, t)|_{L^p} 
\le C t^{-\frac{1}{2}(1-1/p)-1/2},
\ee
\be\label{38}
|{e}_y (y,t)|, |{e}_t (y,t)| \le 
Ct^{-1/2}
\Big(  e^{-\frac{(-y-t)^2} {Ct}} + e^{-\frac{(-y+t)^2} {Ct}}  \Big),
\ee
\be\label{39}
|{e}_{ty} (y,t)| \le C t^{-1} 
\Big(  e^{-\frac{(-y-t)^2} {Ct}} + e^{-\frac{(-y+t)^2} {Ct}}  \Big).
\ee
\el

\begin{proof}
Bound \eqref{35} follows immediately from definition 
\eqref{E}. Given \eqref{38}--\eqref{39},
bounds \eqref{36}--\eqref{37} follow as in
the heat kernel estimates \eqref{heatbds}--\eqref{heatbds2}.
Thus, it remains only to establish \eqref{38}--\eqref{39}.
Differentiating \eqref{E}, we have
$
{e_y}(y,t)=
\Big(\frac{1}{u_+-u_-}\Big)
\Big( \frac { e^{-\frac{(-y-t)^2} {4t}} } {\sqrt{4 \pi t}}
-
\frac { e^{-\frac{(-y+t)^2} {4t}} } {\sqrt{4 \pi t}} \Big) ,
$ 
yielding \eqref{38}(i).
Differentiating \eqref{E} with respect to $t$, we obtain
\begin{equation}\label{te}
\begin{aligned}
{e_t}(y,t)&=
\Big(\frac{-1}{u_+-u_-}\Big)
\Big( \frac { e^{-\frac{(-y-t)^2} {4t}} } {\sqrt{4 \pi t}}
+
\frac { e^{-\frac{(-y+t)^2} {4t}} } {\sqrt{4 \pi t}} \Big) \\&\quad
-
\Big(\frac{t^{-1/2}}{u_+-u_-}\Big)
\Big(  \frac{(-y-t)}{\sqrt{t}}\frac{e^{-\frac{(-y-t)^2} {4t}} } {\sqrt{4 \pi t}}
-
\frac{(-y+t)}{\sqrt{t}}\frac { e^{-\frac{(-y+t)^2} {4t}} } {\sqrt{4 \pi t}} \Big), 
\end{aligned}
\end{equation}
yielding \eqref{38}(ii) immediately for $t\ge 1$.
By the Mean Value Theorem, for $t\le 1$,
\begin{equation}\label{te2}
\begin{aligned}
\Big|  \frac{(-y-t)}{\sqrt{t}}\frac{e^{-\frac{(-y-t)^2} {4t}} } {\sqrt{4 \pi t}}
-
\frac{(-y+t)}{\sqrt{t}}\frac { e^{-\frac{(-y+t)^2} {4t}} } {\sqrt{4 \pi t}} \Big|&=
t \Big|\int_{-1}^1 \partial_z
 \Big(\frac{z}{\sqrt{t}}\frac{e^{-\frac{z^2} {4t}} } {\sqrt{4 \pi t}}\Big)
|_{z=-y+\theta t} \, d\theta\Big|\\
&\le
2Ct
 \Big|\partial_z\Big(\frac{z}{\sqrt{t}}\frac{e^{-\frac{z^2} {4t}} } {\sqrt{4 \pi t}}\Big)|_{z=-y}\Big|\\
&\le
C\Big( e^{-\frac{(-y-t)^2} {Ct}} + e^{-\frac{(-y+t)^2} {Ct}} \Big),  
\end{aligned}
\end{equation}
which, together with \eqref{te}, yields again \eqref{38}(ii).
Estimate \eqref{39} goes similarly.
Note that we have taken crucial account of cancellation in
the small time estimates of $ e_t$, $ e_{ty}$.
\end{proof}

\br\label{2.6}
\textup{
For $t\le 1$, a calculation analogous to 
\eqref{te2} yields
$ | e_y(y,t)|\le Ce^{-\frac{(y+a_-t)^2}{Mt}}, $
and thus $|e(\cdot,s)|_{L^1}\to 0$ as $s\to 0$.
}
\er

\subsection{Bounds on $e^\eps$ }\label{s:ee}

\bl\label{sebds}  For some $C>0$, all $0<\eps\le 1$, and all $t>0$,
\be\label{s35}
|{e^\eps} (\cdot, t)|_{L^\infty}, \le C/\eps ,
\ee
\be\label{s36}
|{e^\eps}_y (\cdot, t)|_{L^p}
\le (C/\eps) t^{-\frac{1}{2}(1-1/p)},
\ee
\be\label{s36b}
|{e^\eps}_t(\cdot, t)|_{L^p} 
\le C t^{-\frac{1}{2}(1-1/p)},
\ee
\be\label{s37}
|{e^\eps}_{ty}(\cdot, t)|_{L^p} 
\le C t^{-\frac{1}{2}(1-1/p)-1/2},
\ee
\be\label{s38}
|{e^\eps}_y (y,t)|\le
(C/\eps)t^{-1/2}
\Big(  e^{-\frac{(-y-t)^2} {Ct}} + e^{-\frac{(-y+t)^2} {Ct}}  \Big),
\ee
\be\label{s38b}
|{e^\eps}_t (y,t)| \le 
Ct^{-1/2}
\Big(  e^{-\frac{(-y-t)^2} {Ct}} + e^{-\frac{(-y+t)^2} {Ct}}  \Big),
\ee
\be\label{s39}
|{e^\eps}_{ty} (y,t)| \le C t^{-1} 
\Big(  e^{-\frac{(-y-t)^2} {Ct}} + e^{-\frac{(-y+t)^2} {Ct}}  \Big).
\ee
\el

\begin{proof}
Bounds \eqref{s35}, \eqref{s36}, and \eqref{s38} follow exactly
as in the $\eps$-independent case.
Bound \eqref{s36b} follows immediately provided that
we can establish \eqref{s38b}, as we now do.
Differentiating \eqref{smallE} with respect to $t$, we obtain
\begin{equation}\label{ste}
\begin{aligned}
{e^\eps_t}(y,t)&=
\Big(\frac{-1}{2}\Big)
\Big( \frac { e^{-\frac{(-y-\eps t)^2} {4t}} } {\sqrt{4 \pi t}}
+
\frac { e^{-\frac{(-y+\eps t)^2} {4t}} } {\sqrt{4 \pi t}} \Big) \\&\quad
-
\Big(\frac{t^{-1/2}}{2\eps}\Big)
\Big(  \frac{(-y-t)}{\sqrt{t}}\frac{e^{-\frac{(-y-\eps t)^2} {4t}} } {\sqrt{4 \pi t}}
-
\frac{(-y+\eps t)}{\sqrt{t}}\frac { e^{-\frac{(-y+t)^2} {4t}} } {\sqrt{4 \pi t}} \Big), 
\end{aligned}
\end{equation}
yielding \eqref{s38b} immediately for $t\ge \eps^{-2}$.
By the Mean Value Theorem, for $t\le \eps^{-2}$,
\begin{equation}\label{ste2}
\begin{aligned}
\Big|  \frac{(-y-\eps t)}{\sqrt{t}}\frac{e^{-\frac{(-y-\eps t)^2} {4t}} } {\sqrt{4 \pi t}}
-
\frac{(-y+\eps t)}{\sqrt{t}}\frac { e^{-\frac{(-y+\eps t)^2} {4t}} } {\sqrt{4 \pi t}} \Big|&=
\eps t \Big|\int_{-1}^1 \partial_z
 \Big(\frac{z}{\sqrt{t}}\frac{e^{-\frac{z^2} {4t}} } {\sqrt{4 \pi t}}\Big)
|_{z=-y+\theta \eps t} \, d\theta\Big|\\
&\le
2C\eps t
 \Big|\partial_z\Big(\frac{z}{\sqrt{t}}\frac{e^{-\frac{z^2} {4t}} } {\sqrt{4 \pi t}}\Big)|_{z=-y}\Big|\\
&\le
C\eps \Big( e^{-\frac{(-y-\eps t)^2} {Ct}} + e^{-\frac{(-y+\eps t)^2} {Ct}} \Big),  
\end{aligned}
\end{equation}
which, together with \eqref{ste}, yields again \eqref{s38b}.
Bounds \eqref{s37} and \eqref{s39} follow similarly.
\end{proof}

%

\subsection{Bounds on $R^\eps$}\label{s:Re}

\bl\label{rt}
For $\cK f:= \int_\RR K(x,y)f(y)\, dy$
and any $1\le p\le \infty$,
\be\label{tKbd}
|\cK f|_{L^p}\le 
\sup_y |K(\cdot,y)|_{L^p} |f|_{L^1},
\ee
\be\label{Kbd}
|\cK|_{L^p\to L^p}\le \max \{
\sup_x |K(x,\cdot)|_{L^1}, \sup_y |K(\cdot,y)|_{L^1} \}
\ee
\el

\begin{proof}
By the Triangle inequality,
$$
\Big| \int_\RR K(\cdot ,y)f(y)dy\Big|_{L^p(x)}\le
 \int_\RR |K(\cdot,y)|_{L^p}|f(y)|dy
\le 
\sup_y |K(\cdot,y)|_{L^p} |f|_{L^1},
$$
establishing \eqref{tKbd}.
This yields also \eqref{Kbd} in case $p=1$.
Likewise,
$$
\Big| \int_\RR K(x,y)f(y)dy\Big|\le
 \int_\RR |K(x,y)|dy |f|_{L^\infty}
\le 
\sup_x |K(x,\cdot)|_{L^1} |f|_{L^\infty},
$$
establishing the \eqref{Kbd} for $p=\infty$.
For general $p$, \eqref{Kbd}
then follows by the 
Riesz--Thorin Interpolation Theorem.
\end{proof}

\bl\label{Rbds}  For some $C>0$, all $0< \eps\le 1$, and all $t>0$,
\be\label{sRbdeq}
\sup_y |{R^\eps} (\cdot, t;y)|_{L^p(x)} ,
\quad
\sup_x |{R^\eps} (x, t; \cdot)|_{L^p(y)} \le Ct^{-\frac{1}{2}(1-1/p)},
\ee
\be\label{sRbdeqy}
\sup_y |{R^\eps_y} (\cdot, t;y)|_{L^p(x)} ,
\quad
\sup_x |{R^\eps_y} (x, t; \cdot)|_{L^p(y)} \le 
Ct^{-\frac{1}{2}(1-1/p)-\frac{1}{2}}.
\ee
\el

\begin{proof}
From $(\bar u^\eps)'\sim \eps^2 e^{-\theta \eps |x|}$, we obtain 
\ba\label{sRcalc}
R^\eps (x,t;y) & =
(1/2\eps) x (\bar u^\eps)'(x) \int_0^1
\Big(\frac {e^{-\frac{(\theta x-y-\eps t)^2}{4t}} }{\sqrt{4 \pi t}} 
-\frac {e^{-\frac{(\theta x-y+\eps t)^2}{4t}} }{\sqrt{4 \pi t}} \Big)
d\theta\\
&\le
Ce^{-\theta \eps |x|} \int_0^1
\Big(\frac {e^{-\frac{(\theta x-y-\eps t)^2}{4t}} }{\sqrt{4 \pi t}} 
+\frac {e^{-\frac{(\theta x-y+\eps t)^2}{4t}} }{\sqrt{4 \pi t}} \Big)
d\theta,
\ea
from which we obtain immediately
$|R^\eps|_{L^\infty}\le Ct^{-1/2}$, and, bounding 
$Ce^{-\theta \eps |x|}$ by $C$,
$$
\sup_x |R^\eps|_{L^p(y)}\le Ct^{-\frac{1}{2}(1-1/p)}
$$
 for any $p$.

Bounding the integral on the righthand side
 by $C_1t^{-1/2}$ and the $L^1(x)$ norm
of $Ce^{-\theta \eps |x|}$ by $C_2/ \eps$,
we find
$\sup_y |R^\eps|_{L^1(x)}\le C_2 t^{-1/2}/\eps\le C$
for $t\ge \eps^{-2}$.
For $t\le \eps^{-2}$, on the other hand, we may estimate
the integral (the middle displayed term in the first equality) instead,
using the Mean Value Theorem, as
$$
\begin{aligned}
\int_0^1
\Big(\frac {e^{-\frac{(\theta x-y-\eps t)^2}{4t}} }{\sqrt{4 \pi t}} 
-\frac {e^{-\frac{(\theta x-y+\eps t)^2}{4t}} }{\sqrt{4 \pi t}} \Big)\,d\theta
&\le 
\int_0^1
(2\eps t)
\partial_z \Big(\frac {e^{-\frac{(\theta x-y-z)^2}{4t}} }{\sqrt{4 \pi t}} 
\Big)|_{z=z_*\in [-\eps t, \eps t]}\, d\theta\\
&\le
(2\eps t) 
\int_0^1 Ct^{-1} \, d\theta \le C\eps,
\end{aligned}
$$
to again obtain $\sup_y |R^\eps|_{L^1(x)}\le C_2 \eps/\eps\le C$.
The bounds on
$\sup_y |R^\eps|_{L^p(x)}$ then follow by H\"older interpolation
between the $L^1$ and $L^\infty$ bounds, verifying \eqref{sRbdeq}
Similar computations yield \eqref{sRbdeqy}.
\end{proof}

\subsection{Euler--Lagrange equations for least squares}\label{s:el}

Setting $E(\alpha):= \frac{1}{2}|u|_{L^2}^2
= \frac{1}{2}| \tilde u(\cdot +\alpha, t)-\bar u(\cdot ) |_{L^2}^2 $
and differentiating, we have
$$
\begin{aligned}
\frac{dE}{d\alpha}&=
\langle \tilde u(\cdot +\alpha, t)-\bar u(\cdot),
\tilde u'(\cdot + \alpha, t) \rangle_{L^2}
= \langle u, \bar u' + u'\rangle_{L^2}
=
 \langle u, \bar u'\rangle_{L^2},
\end{aligned}
$$
where, in the final equality, we have
used $\langle u, u'\rangle_{L^2}=
\int_{-\infty}^{+\infty} (u^2/2)'(x) dx=0 $
for $u\in H^1$.

\end{document}